\documentclass[12pt, reqno]{amsart}
\usepackage{amssymb}
\usepackage{amsfonts}
\usepackage{amssymb,amsmath,amscd}
\usepackage[colorlinks=true]{hyperref}
\usepackage{mathrsfs}
\newtheorem{Theorem}{ Theorem}

\newtheorem{Le}[Theorem]{Lemma}
\newtheorem{Corollary}[Theorem]{Corollary}

\newtheorem{Proposition}[Theorem]{Proposition}

\theoremstyle{remark}
\newtheorem{Remark}{Remark}

\newcommand{\f}{\frac}

\newcommand{\al}{\alpha}
\newcommand{\de}{\delta}
\newcommand{\va}{\varepsilon}
\newcommand{\mn}{\mathbb{N}}
\newcommand{\ta}{\theta}
\newcommand{\sst}{\subset}
\newcommand{\sls}[1]{\sum_{#1}}
\newcommand{\eq}{\equiv}
\newcommand{\lr}[1]{\left({#1}\right)}
\newcommand{\rd}{\rho_d}
\newcommand{\rdo}{\rho_{d_1}}
\newcommand{\rdt}{\rho_{d_2}}
\newcommand{\mj}{\mathcal{J}^{-1}}
\newcommand{\ra}{\rho^+}
\newcommand{\mr}{\mathbb{R}}
\newcommand{\mt}{\mathbb{T}}
\newcommand{\mz}{\mathbb{Z}}

\newcommand{\mm}{\mathcal{M}}
\newcommand{\mk}{\mathfrak{m}}
\newcommand{\sw}{\sigma_W(b)}
\newcommand{\ew}[1]{e_W({#1})}
\newcommand{\xw}{\frac{x}{dW}}
\newcommand{\xyw}{\frac{x+y}{dW}}
\newcommand{\vq}{V_q(a[d_1,d_2]^k,b,[d_1,d_2],0)}
\newcommand{\me}{\mathbb{E}}

\newcommand{\my}{\mathcal{Y}}
\newcommand{\RNum}[1]{\uppercase\expandafter{\romannumeral #1\relax}}
\renewcommand{\bar}{\overline}

\newcommand{\md}[1]{\,(\textrm{mod}\,{#1})}

 \newcommand{\set}[1]{\left\{#1\right\}}
\newcommand{\bigset}[1]{\bigl\{ #1 \bigr\}}
\newcommand{\Bigset}[1]{\Bigl\{ #1 \Bigr\}}

\newcommand{\bigabs}[1]{\bigl| #1 \bigr|}
\newcommand{\Bigabs}[1]{\Bigl| #1 \Bigr|}

\newcommand{\Bigfloor}[1]{\Bigl\lfloor #1 \Bigr\rfloor}
\newcommand{\floor}[1]{\left\lfloor #1 \right\rfloor}

\newcommand{\bigbrac}[1]{\bigl( #1 \bigr)}
\newcommand{\Bigbrac}[1]{\Bigl( #1 \Bigr)}
\newcommand{\biggbrac}[1]{\biggl( #1 \biggr)}

\newcommand{\norm}[1]{\left\| #1\right\|}

\newcommand{\E}{\mathbb{E}}
\newcommand{\Z}{\mathbb{Z}}
\newcommand{\T}{\mathbb{T}}
\newcommand{\C}{\mathbb{C}}

\newcommand{\R}{\mathbb{R}}

\renewcommand{\hat}[1]{\widehat{#1}}
\renewcommand{\epsilon}{\varepsilon}

\begin{document}

\title{
 Waring-Goldbach problem in short intervals
}

\author{Mengdi Wang}


\address{School of Mathematics\\ Shandong University\\Jinan  250100 \\ China}

\curraddr{Department of Mathematics\\ KTH\\ Stockholm  100 40\\ Sweden}

\email{mengdiw@kth.se}

\maketitle

\begin{abstract}
Let $k\geq2$ and $s$ be positive integers. Let $\ta\in(0,1)$ be a real number. In this paper, we establish that if $s>k(k+1)$ and $\ta>0.55$, then every sufficiently large natural number $n$, subjects to certain congruence conditions, can be written as 
$$
n=p_1^k+\cdots+p_s^k,
$$
where $p_i(1\leq i\leq s)$ are primes in the interval $((\f{n}{s})^{\f{1}{k}}-n^{\f{\ta}{k}},(\f{n}{s})^{\f{1}{k}}+n^{\f{\ta}{k}}]$.  The second result of this paper is to show that if $s>\f{k(k+1)}{2}$ and $\ta>0.55$, then almost all  integers $n$, subject to certain congruence conditions, have above representation.
\end{abstract}

 \section{Introduction}
Suppose that $k\geq2$ is a positive integer and $p$ is a prime. Let $\tau=\tau(k,p)$ be the integer such that $p^\tau\|k$, which means that $p^{\tau}|k$ but $p^{\tau+1}\nmid k$. Define $\gamma=\gamma(k,p)$ by taking
\begin{align}\label{ga}
\gamma=
\begin{cases}
\tau+2\qquad &\text{if\ }p=2\ \text{and}\ \tau>0;\\
\tau+1	 &\text{otherwise.}
\end{cases}
\end{align}
We now put
\begin{align}\label{rk}
R_k=\prod_{(p-1)|k}p^{\gamma}.	
\end{align}
Waring-Goldbach problem asks whether every sufficiently large number $n$ with $n\eq s\pmod{R_k}$ can be expressed as a sum of $s$ terms of $k$-th powers of primes, where $s$ is a positive number depending on $k$. The first result was obtained by Hua \cite{Hua38}, who showed that when $s\geq2^{k}+1$ every sufficiently large natural number $n$ with $n\eq s\pmod{R_k}$ can be represented as
$$
n=p_1^k+\cdots+p_s^k,
$$
where $p_1,\cdots,p_s$ are prime numbers. Subsequent works are focus on reducing the value of $s$. The latest development for general $k$ is that in 2017, Kumchev and Wooley \cite{KW} showed that above statement is true for large values of $k$ and $s>(4k-2)\log k-(2\log 2-1)k-3$.

One interesting generalization of Waring-Goldbach problem is to restrict the prime variables into a short interval. Set $0<\theta_{k,s}<1$ such that for each $\ta>\ta_{k,s}$ the equation system 
\begin{align}\label{eq}
\begin{cases}
n=p_1^k+\cdots+p_s^k;\\
	|\lr{\f{n}{s}}^{1/k}-p_i|\leq\lr{\f{n}{s}}^{\f{\ta}{k}}\ (1\leq i\leq s)
\end{cases}
\end{align}
is solvable for all sufficiently large $n$ with $n\eq s\pmod{R_k}$. In 2015, Wei and Wooley \cite{WW} proved that when $s>\max\{6,2k(k-1)\}$,
$$
\ta_{k,s}\leq	
\begin{cases}
\f{19}{24}\qquad &\text{if\ }k=2;\\
\f{4}{5}&\text{if\ }k=3;\\
\f{5}{6}&\text{if\ }k\geq4.
\end{cases}
$$
This result has been improved by Huang \cite{Huang}, Kumchev and Liu \cite{KL}, Matom\"{a}ki and Shao \cite{MS} successively. Matom\"{a}ki and Shao showed that when $k\geq2$ and $s>k(k+1)+2$ one could get $\ta_{k,s}\leq\f{2}{3}$. So far as we are concerned the improvements to $\ta_{k,s}$ always rely on two aspects: either estimating
\begin{align}\label{es}
\sls{X<x\leq X+H}\Lambda(x)e(x^k\al)	
\end{align}
for a smaller $H$ directly, such as \cite{MS}, where $\Lambda$ is the von Mangoldt function; or considering using (\ref{es}) in a more efficient way, such as \cite{Wang}, by making use of an argument in \cite{Zhao}, also in \cite{KZ}. Recently, Salmensuu \cite{Sal} applied the transference principle to this problem to prove that when $s>k(k+1)$,
$$
\ta_{k,s}\leq	
\begin{cases}
0.644\qquad &\text{if\ }k=2;\\
0.578&\text{if\ }k=3;\\
0.55&\text{if\ }k\geq4.
\end{cases}
$$

Green \cite{Gr} proposed the transference principle to handle translation-invariant additive problems. Matom\"{a}ki, Maynard, Shao \cite{MMS} and Salmensuu \cite{Sal} extended the transference principle to general (non translation-invariant) additive problems. Compared with previous results in Waring-Goldbach problem with almost equal summands, the application of transference principle cannot give the asymptotic formula for counting the number of solutions to (\ref{eq}). However, on the other hand, the transference principle would reduce the dependence of exponential sum (\ref{es}). In practice, we just need to show that for each frequency $0\neq\alpha\in\mathbb R/\mathbb Z$,
\begin{align}\label{eg}
\sls{X<x\leq X+H}\nu(x)e(x^k\al)=o(H),
\end{align}
for a suitable sieve majorant $\nu$. The freedom of the choice of majorant makes it possible to avoid any complicated computations in major arcs such as the strategy of enlarging major arcs. We would lead readers to a survey \cite{Pre} written by Prendiville for a detailed introduction to various version of the transference principle.  To our knowledge, in order to estimate (\ref{es}), the classical method of applying Vaughan's identity reduces the task to dealing with both type $\RNum{1}$ sum and type $\RNum{2}$ sum. But, it is not hard to find in the following (indeed in Lemma 5) that we just need to consider type $\RNum{1}$  sum in estimating (\ref{eg}). This also simplifies the exponential sum estimates in minor arcs. Contrasted with \cite{Sal}, we obtain a smaller $H$ in the price of a worse upper bound of (\ref{eg}) when $\alpha$ is in minor arcs. However, we can afford this because the dominant upper bound of (\ref{eg}) is determined by the major arcs.

\begin{Theorem}
Suppose that $2\leq k<s$ are positive integers. Suppose that $\va>0$  and $\ta\in(\f{1}{2},1)$ are two parameters. Let $\al^-$ be a positive number such that for every large enough $x$, every interval $I\sst[x,x+x^{\ta+\va}]$ with length $|I|\geq x^{\ta-\va}$, and every pair of integrs $c,d\in\mn$ with $(c,d)=1$ and $d\leq\log x$ we always have
\begin{align}\label{ww}
\sls{\mbox{\tiny$\begin{array}{c}
p\in I\\
p \text{\ is a prime}\\
p\equiv c\pmod d\end{array}$}}	1\geq \f{\al^-|I|}{\phi(d)\log x}.
\end{align}
Then when 
$$
s>\max\{\f{2}{\al^-\ta},k(k+1)\}
$$
for every sufficiently large integer $n\eq s\pmod{R_k}$, where $R_k$ is defined in (\ref{rk}), there exist primes $p_1,\cdots,p_s\in[\lr{\f{n}{s}}^{1/k}-n^{\f{\ta}{k}},\lr{\f{n}{s}}^{1/k}+n^{\f{\ta}{k}}]$ such that
$$
n=p_1^k+\cdots+p_s^k.
$$

\end{Theorem}

Thanks to \cite[Theorem 10.3]{Har} and \cite[Theorem 10.8]{Har}, we have the following corollaries.
\begin{Corollary}
\begin{enumerate}
\item When $k>7, s>k(k+1)$ and $\theta>0.525$, for all sufficiently large $n\eq s\pmod{R_k}$, there are primes $p_1,\dots,p_s\in[\lr{\f{n}{s}}^{1/k}-n^{\f{\ta}{k}},\lr{\f{n}{s}}^{1/k}+n^{\f{\ta}{k}}]$ such that
\[
n=p_1^k+\cdots+p_s^k; 
\]
\item When $2\leq k\leq7, s>k(k+1)$ and $\theta>0.55$, for all sufficiently large $n\eq s\pmod{R_k}$, there are primes $p_1,\dots,p_s\in[\lr{\f{n}{s}}^{1/k}-n^{\f{\ta}{k}},\lr{\f{n}{s}}^{1/k}+n^{\f{\ta}{k}}]$ such that
\[
n=p_1^k+\cdots+p_s^k;
\]
\item When $2\leq k\leq7, s>43$ and $\theta>0.525$, for all sufficiently large $n\eq s\pmod{R_k}$, there are primes $p_1,\dots,p_s\in[\lr{\f{n}{s}}^{1/k}-n^{\f{\ta}{k}},\lr{\f{n}{s}}^{1/k}+n^{\f{\ta}{k}}]$ such that
\[
n=p_1^k+\cdots+p_s^k.
\] 	
\end{enumerate}
\end{Corollary}

It is worth noting that our Theorem 1 improves \cite[Theorem 1]{Sal} in both the case $\theta_{k,s}\leq 0.55$ when $k=2,3$ and extreme case $\theta_{k,s}\leq 0.525$ when $2\leq k\leq22$.

The current knowledge can not handle the equations with few variables, such as the Goldbach conjecture. However, in 1938,  Chudakoff \cite{Chu}, van der Corput \cite{van} and Esterman \cite{Est} showed that all large integers $n\leq N$ out of an exceptional set with size $o(N)$ are sums of two primes. This encourages us to study the almost all version of the Waring-Goldbach problem with almost equal variables. In the traditional circle method, this is achieved by using an average idea, precisely, using Bessel's inequality (see \cite[Section 9]{WW} as an example).  And this allows us to reduce roughly half as many variables in the equation. Here we handle the exceptional set result by establishing an almost all version of the transference principle.

\begin{Theorem}
Suppose that $2\leq k<s$ are positive integers and $M>0$ is a large number. Suppose that $\va>0$  and $\ta\in(\f{1}{2},1)$ are two parameters.  Let $\al^-$ be the parameter defined as in Theorem 1. Then when $s>\max\{\f{2}{\al^-\ta},\f{k(k+1)}{2}\}$, for all but $O\bigbrac{(\log_4 M)^{-\frac{1}{12s+6}} M}$\footnote{The subscript indicates the number of iterated logarithms.} of $n\leq M$ subjects to 	$n\eq s\pmod{R_k}$ (besides, in the case $k=3$ and $s=7$ we also need $n\not\equiv0\pmod{9}$), one can find primes $p_1,\cdots,p_s$ in the interval $|\lr{\f{n}{s}}^{1/k}-p_i|\leq n^{\f{\ta}{k}}(1\leq i\leq s)$ such that
$$
n=p_1^k+\cdots+p_s^k.
$$
\end{Theorem}

We  point out that the additional congruence condition $n\not\equiv0\pmod{9}$  is to ensure $p_1^3+\cdots+p_7^3=n$  solvable. Now, let $M$ be a sufficiently large positive number. Let $\ta_{k,s}'$ be the smallest $\ta$ such that the following equation is solvable for almost all $n\in(M-(\f{M}{s})^{\f{k-1+\ta}{k}},M]$
$$
\begin{cases}
n=p_1^k+\cdots+p_s^k;\\
	|\lr{\f{M}{s}}^{1/k}-p_i|\leq M^{\f{\ta}{k}}\ (1\leq i\leq s).
\end{cases}
$$
One then can deduce from Theorem 2 and \cite[Theorem 10.3]{Har} that $\ta'_{k,\f{k(k+1)}{2}+1}\leq0.55+\va$ when $k\geq2$. This improves \cite[Theorem 2]{KL} which shows  that $\ta'_{k,\f{k(k+1)}{2}+1}\leq\f{31}{40}+\va$.

\subsection*{Acknowledgements}

The author would like to express her gratitude to her advisor Xuancheng Shao for many helpful discussions, suggestions and the modification of the draft; to her advisor Lilu Zhao for checking the draft and giving useful suggestions; to  Juho Salmensuu for his detailed and helpful comments.
The author is also obliged to the anonymous referee for the helpful, detailed and instructive comments, especially simplified the proof of Theorem 3, and helped to construct a quantitative Proposition \ref{almost}.
 The author would like to thank the financial support from China Scholarship Council for supporting her stay in the US, and also thank and the Department of Mathematics at the University of Kentucky for the hospitality and the excellent conditions.

\bigskip

\section{Notations}

We would like to introduce  some basic notation common to the whole paper here. Let $\va>0$ be a small number  and it  is  allowed to change at different occurrences. With or without subscript, letters $p$ and $p_i(i=1,2,\cdots)$ denote prime numbers. If $X$ is a positive integer, we would write $[X]$ for the discrete interval $\{1,2,\cdots,X\}$. We would also abbreviate $e^{2\pi ix}$ to $e(x)$.

Assume that $f:\mr\rightarrow\mathbb{C}$ and $g:\mr\rightarrow\mr_{\geq0}$ are two functions, we make the following notations:

\begin{itemize} 
\item $f\ll g$ means that there exists some constant $C>0$ and a real number $x_0$ such that for all $x\geq x_0$ we have $|f(x)|\leq Cg(x)$;
\item $f=O(g)$ is the same as $f\ll g$;
\item $f\asymp g$ means that $f\ll g$ and $g\ll f$ ($f$ should be  real-valued and positive);
\item $f=o(g)$ if $g\neq0$ and $\lim_{x\rightarrow\infty}\f{f(x)}{g(x)}=0$.	
\end{itemize}

We'll use counting measure on $\Z$, so for a function $f:\Z\to\C$ its $L^p$-norm is defined to be
\[
\norm{f}_p^p=\sum_x|f(x)|^p,
\]
and $L^\infty$-norm is defined to be
\[
\norm{f}_\infty=\sup_x|f(x)|.
\]
Besides, we'll use Haar probability measure on $\T$, so for a function $F:\T\to\C$ define its $L^p$-norm as
\[
\norm{F}_p^p=\int_\T|F(\alpha)|^p\,\mathrm d\alpha.
\]

If $f:A\rightarrow\mathbb{C}$ is a function and $B$ is a non-empty finite subset of $A$, we write
$$
\E_{x\in B}f(x)=\frac{1}{|B|}\sum_{x\in B}f(x)
$$
as the average of $f$ on $B$. We would also abbreviate $\E_{x\in A}f(x)$ to $\E_A(f)$ if the supported set $A$ is finite and no confusion caused.

We'll use Fourier analysis on $\Z$ with its dual $\T$. Let $f:\Z\to\C$ be a function, define its Fourier transform by setting
\[
\hat{f}(\alpha)=\sum_xf(x)e(x\alpha),
\]
where $\alpha\in\T$. For functions $f,g:\Z\to\C$ define the convolution of $f$ and $g$ as
\[
f*g(x)=\sum_yf(y)\bar{g}(x-y).
\]

\bigskip

\section{Proof of Theorem 1}

Suppose that $n$ is a sufficiently large integer and subjects to the congruence condition $n\eq s\pmod{R_k}$, where $R_k$ is defined in (\ref{rk}). Write $\tilde x=\lceil{n/s\rceil}^{\frac{1}{k}}$ and $\tilde y=\tilde x^\theta$, where $0<\theta<1$ is a parameter. We care about  whether the below equation has prime solutions $p_1\dots,p_s$ in the interval $[\tilde x-\tilde y,\tilde x-\tilde y]$
\[
n=p_1^k+\cdots+p_s^k.
\]
Clearly, since we expect each $p_i^k$ surrounds $\frac{n}{s}$, the parameter $\theta$ is to capture how short the interval would be such that above equation is solvable. In order to apply the transference principle, we are going to restrict those primes $p_i$ subject to a congruent equation, named $p_i^k\equiv b\pmod{W}$, where $1\leq b\leq W$ are numbers defined as follows, 
\[
w=\log\log\log n, \quad W=2k^2\lceil{\kappa^{-1}\rceil}!^2\prod_{p\leq w}p,\quad (b,W)=1
\]
where $\kappa\in(0,1)$ is a sufficiently small number. We call the procedure to restrict the $k$-th power of primes $p^k$  into a reduced residue class modulus $W$  as the $W$-trick. The next step, let
\begin{align}\label{yy}
m&=\Bigl\lfloor{\frac{(\tilde x-\tilde y)^k}{W}\Bigr\rfloor}, \quad   N=\Bigl\lceil{\frac{(\tilde x+\tilde y)^k-(\tilde x-\tilde y)^k}{W}\Bigr\rceil},\nonumber\\
X&=Wm+b,  \quad\qquad Y=WN.
\end{align}
It is not hard to find that
\begin{align}\label{xy}
Y\asymp\tilde x^{k-1}\tilde y\asymp X^{1-\f{1}{k}+\f{\ta}{k}},	
\end{align}
and
\begin{align}\label{w}
	W\ll\log\log X. 
\end{align}

 In addition,  We would construct the standard Selberg's sieve majorant. Assume that $0<\delta<1$ is a parameter, following the notations in \cite[p.15]{Shao}, set
  \begin{align}\label{z}
	z=X^{\de/2},\quad D=z^2\text{ and }P=\prod_{w\leq p\leq z}p.
	\end{align}
And let $\rho_d$ be weights which are supported on $d|P$ with $d<z$, and satisfy  $|\rd|\leq1$  for all $d$ and $\rho_1=1$. Besides, set $y_d=\mu(d)\phi(d)\sum_{d|m}\frac{\rho_m}{m}$ which satisfies $y_d=\mathcal J^{-1}$, where
\[
\mathcal J=\sum_{d<z\atop (d,W)=1}\frac{1}{\phi(d)},
\]
and define
\begin{align}\label{r}
\ra(n)=\lr{\sls{d|(n,P)}\rd}^2,
\end{align}
and 
\begin{align}\label{1}
\al^+=\f{\phi(W)}{kW}\log X\sls{d_1,d_2|P}\f{\rdo\rdt}{[d_1,d_2]}.
\end{align}
Then one can deduce from \cite[(A.1)]{Shao}, \cite[Lemma A.3]{Shao} and (\ref{z})  that
\begin{align}\label{al}
\al^+=\f{1}{k}\f{\log X}{\log z}+O_W(\f{1}{\log X})=\f{2}{k\de}+O_W(\f{1}{\log X}).	
\end{align}
Now let $f_b,\nu_b:[N]\rightarrow\mr_{\geq0}$ be functions defined as \cite[(26)-(28)]{Sal}, that is,
\begin{align}\label{f}
&f_b(n)=
\begin{cases}
\f{\phi(W)}{\al^+W\sigma_W(b)}X^{1-\f{1}{k}}\log X	\qquad &\text{if\ }W(m+n)+b=p^k \text{\ for some }p;\\
0 &\text{otherwise,}
\end{cases}\\
\label{vb}
&\nu_b(n)=\begin{cases}
\f{\phi(W)}{\al^+W\sigma_W(b)}X^{1-\f{1}{k}}\log X	\ra(t)\qquad &\text{if\ }W(m+n)+b=t^k ;\\
0 &\text{otherwise,}
\end{cases} 
\end{align}
and 
\begin{align}\label{si}
\sigma_W(b)=\#\{z\in[W]:z^k\eq b\pmod{W}\}.	
\end{align}

Here the function $f_b$ is a weighted version of indicator of $p^k$. We would like to do some explaining to make it a little bit easy to follow. Combining the restrictions of primes, that are $(\tilde x-\tilde y)^k\leq p^k\leq (\tilde x+\tilde y)^k$ and $p^k\equiv b\md{W}$, we may rewrite $p^k$ as $p^k=W(m+n)+b$, where $Wm$ is the left endpoint $(\tilde x-\tilde y)^k$. The factor $\log X$ in $f_b$ is a weight for primes; $X^{1-\frac{1}{k}}$ is a linearity factor, just like the setting of \cite{BP} and \cite{Chow}; and $\f{\phi(W)}{\al^+W\sigma_W(b)}$ is a normalized factor to ensure the majorant $\nu_b$ has $L^1$-norm $N$.

\begin{Le}[Pseudorandomness] 
Let $\al\in\mt$, $1/2<\ta<1$ and $2\leq k\in\mathbb N$. Let $\nu_b:[N]\rightarrow\mr_{\geq0}$ be the function defined in (\ref{vb}).
When $\de<\f{\ta}{k}$ (recall (\ref{z})) we can get
$$
|\widehat\nu_b(\alpha)-\widehat{1_{[N]}}(\alpha)|\ll N(\log\log\log  N)^{-\frac{1}{2k}}.
$$
\end{Le}

We now apply transference principle, see \cite[Lemma 6]{Sal}, together with  \cite[Lemma 7, 9]{Sal} and above Lemma 4 to prove the first main theorem.

\vspace{3mm}

\noindent\emph{Proof of Theorem 1.}

\vspace{2mm}

Let $s\geq3k$, suppose that  $n$ is a sufficiently large integer and subjects to $n\eq s\pmod{R_k}$. By taking $q=W$ in \cite[Lemma 10]{Sal} and then apply this lemma, one can find a group of integers $1\leq b_1,\cdots, b_s\leq W$ with $(b_1\cdots b_s,W)=1$ such that
$$
n\eq b_1+\cdots+b_s\pmod{W}.
$$
One can use these integers $b_1,\cdots, b_s$ to define functions $f_{b_1},\dots,f_{b_s}$, and majorants $\nu_{b_1},\dots,\nu_{b_s}$ in the manner of (\ref{f}) and (\ref{vb}). \cite[Lemma 6 -Lemma 7, Lemma 9]{Sal} together with Lemma 4 reveal that when $s>\max\{\f{2}{\al^-\ta},k(k+1)\}$, for any $n_0\in[\f{N}{2},N]$ we have
\[
f_{b_1}*\dots*f_{b_s}(n_0)\gg N^{s-1}.
\]
That is, from the definition of functions $f_{b_i}$, when $s>\max\{\f{2}{\al^-\ta},k(k+1)\}$, for any $n_0\in[\f{N}{2},N]$ there exist positive integers $n_1,\cdots,n_s$ and primes $p_1\cdots,p_s$ such that
$$
\begin{cases}
n_0=n_1+\cdots+n_s\\
W(m+n_i)+b_i=p_i^k(1\leq i\leq s).	
\end{cases}
$$

On the other hand, it can be seen that $\frac{n-b_1-\dots-b_s-smW}{W}\in[\frac{N}{2},N]$. Thus, employing above conclusion to $\frac{n-b_1-\dots-b_s-smW}{W}$ shows that when $s>\max\{\f{2}{\al^-\ta},k(k+1)\}$,  we can find a group of integers $n_1,\cdots,n_s$  and a group of primes $p_1\cdots,p_s$ such that
\[
\begin{cases}
n=W(sm+n_1+\dots+n_s)+b_1+\dots+b_s\\
W(m+n_i)+b_i=p_i^k(1\leq i\leq s).	
\end{cases}
\]
Then the theorem follows by verifying that primes $p_1\cdots,p_s$ are in the interval $[\tilde x-\tilde y,\tilde x+\tilde y]$.

\qed

\vspace{3mm}

In order to prove  Lemma 4 we shall use the Hardy-Littlewood method. We would discuss both cases divided according to whether $\al$ nears a rational number with small denominator or not. Let
\begin{align}\label{q}
Q=(\log X)^{c_kA},
\end{align}
where $c_k$ is a positive number only depending on the degree $k$, $A\geq20$ is a sufficiently large positive number. We now define the major arcs as
\[
\mm =\bigcup_{q\leq Q}\bigcup_{1\leq a\leq q\atop(a,q)=1}\mm(q,a)\quad\text{ and  }\mm(q,a)=\{\al:|\al-\f{a}{q}|\leq\f{Q}{qY}\},
\]
where $Y$ is defined in (\ref{yy}), and  minor arcs are the complement of major arcs $\mk =\mt\backslash\mm.$

Recalling the definition of $\nu_b$ in (\ref{vb}) and recalling (\ref{r}), one can deduce from Fourier analysis that for any $\al\in\mt$,
\begin{align}\label{v}
	\widehat\nu_b(\al)&=\sls{n}\nu_b(n)e(n\al)\nonumber\\
	&=e\lr{-(\f{b}{W}+m)\al}\f{\phi(W)}{\al^+W \sw} X^{1-\f{1}{k}}\log XE_b(\al),
\end{align}
where
\begin{align}\label{e}
E_b(\al)=\sls{X<t^k\leq X+Y\atop t^k\equiv b\md{W}}\lr{\sls{d|(P,t)}\rd}^2\ew{t^k\al},
\end{align}
and $\ew{\cdot}=e(\f{\cdot}{W})$.

We are going to apply \cite[Proposition 2.1]{MS} (a minor arc type $\RNum{1}$ estimate) to deal with the exponential sum $E_b(\al)$ when $\al$ belongs to minor arcs in Lemma 5; and in Lemma 7, we would do the major arcs analysis in a similar manner of \cite[Section 7]{Sal}.

\begin{Le}[Minor arc estimates]
Let $\ta\in(\f{1}{2},1)$ and $k\geq2$, when $\de<\f{\ta}{k}$ we have
$$
\sup_{\al\in\mk}|\widehat\nu_b(\al)-\widehat{1_{[N]}}(\al)|\ll\frac{N}{\log N}.
$$
\begin{proof}
Clearly,
\[
|\widehat1_{[N]}(\alpha)|=|\sum_{n\leq N}e(n\alpha)|\ll\norm{\alpha}^{-1}.
\]
The making use of Dirichlet approximate theorem shows that for every $\alpha\in\T$ there is an integer $1\leq q\leq Q$ such that $\norm{q\alpha}\leq\frac{1}{Q}$. As a consequence, when $\alpha\in\mk$, there is $1\leq q\leq Q$ such that $\frac{Q}{Y}\leq\norm{q\alpha}\leq\frac{1}{Q}$. Thus,
\[
|\widehat1_{[N]}(\alpha)|\ll\norm{\alpha}^{-1}\ll\norm{q\alpha}^{-1}\ll\frac{N}{\log N}.
\]

Therefore, we just need to show that $\widehat\nu_b(\al)=O(\frac{N}{\log N})$ whenever $\al$ belongs to minor arcs. By (\ref{yy}), (\ref{xy}), (\ref{al}) and (\ref{v}), it suffices to prove that when $\al\in\mk$
$$
E_b(\al)\ll X^{\f{\ta}{k}}(\log X)^{-5}.
$$
With the help of \cite[Proposition 2.1]{MS}, we will prove this result by contradiction. Let $\omega=(\log X)^{-A}$, where $A\geq20$ is a large positive  number. And assume that
\begin{align}\label{3}
\sw\f{\omega X^{\f{\ta}{k}}}{W}\lr{\log X}^{14}	\leq  |E_b(\al)|.
\end{align}
We aim to show that (\ref{3}) would never happen. To begin with, from the definition of $E_b(\alpha)$ in (\ref{e}), one can expand the square in $E_b(\alpha)$ to get that
$$
E_b(\al)=\sls{X<t^k\leq X+Y\atop t^k\equiv b\md{W}}\sls{d_1,d_2|(P,t)}\rdo\rdt\ew{t^k\al}.
$$
Exchanging the order of summation gives that
$$
E_b(\al)=\sls{d_1,d_2|P}\rdo\rdt\sum_{\substack{
 X<t^k\leq X+Y\\
 [d_1,d_2]|t\\
 t^k\equiv b\md{W} }}
 \ew{t^k\al}.
$$
For simplicity of writing, now let $x=\lfloor{X^{\f{1}{k}}\rfloor}, y=(X+Y)^{\f{1}{k}}-\lfloor{X^{\f{1}{k}}\rfloor}$. It is not hard to find from (\ref{xy}) that
$y\asymp x^\ta$. Let $d=[d_1,d_2]$ be the least common multiple of $d_1,d_2$, one can see from $d_1,d_2|P$ that $(d,W)=1$; besides we have $d|P$ dues to $d_1|P$, $d_2|P$ and $P$ is square-free. Taking note that $\#\{(d_1,d_2):\text{lcm}[d_1,d_2]=d\}<(\tau(d))^2<\tau_3(d)$, where $\tau_l$ is the $l$-fold divisor function, as well as noting that $|\rho_d|\leq 1$ and $\rho_d$ is  supported in the range $d\leq D$, one has
\begin{align*}
|E_b(\al)|&\leq\sls{d|P\atop d\leq D}\tau_3(d)\sls{z\in[W]\atop z^k\equiv b\md{W}}\bigg|\sls{x<dl\leq x+y\atop dl\equiv z\md{W}}\ew{d^kl^k\al}\bigg|\\
&=\sls{z\in[W]\atop z^k\equiv b\md{W}}\sls{d|P\atop d\leq D}\tau_3(d)\bigg|\sls{\f{x}{d}<l\leq\f{x+y}{d}\atop l\equiv z\bar d\md{W}}\ew{d^kl^k\al}\bigg|,
\end{align*}
where we have written $\bar d$ as the inverse of $d\pmod{W}$. We now assume that $l=z\bar d+uW$, and then divide the summation over $d$ into dyadic intervals to get
$$
|E_b(\al)|\ll \log x\sls{z\in[W]\atop z^k\equiv b\md{W}}\sls{d\sim M}\tau_3(d)\bigg|\sls{\xw-\f{z\bar d}{W}<u\leq\xyw-\f{z\bar d}{W}}\ew{(zd\bar d+Wdu)^k\al}\bigg|,
$$
where $M\leq D$. Combining above inequality with assumption (\ref{3}), one can see from pigeonhole principle that there is some $z\in[W]$ such that
$$
\f{\omega y}{W}(\log x)^{13}\ll\Bigg|\sls{d\sim M}\tau_3(d)c_d\sls{\xw-\f{z\bar d}{W}<u\leq\xyw-\f{z\bar d}{W}}\ew{(zd\bar d+Wdu)^k\al}\Bigg|,
$$
where $c_d\in\mathbb C$ and $|c_d|=1$.
Write $g(du)$ as above phase function, i.e. $g(du)=\frac{(zd\bar d+Wdu)^k\al}{W}$. Plainly, we have $g(du)=\frac{\al\big((zd\bar d+Wdu-x)+x\big)^k}{W}$. It then follows from binomial theorem that
\begin{align*}
g(du)&=\sls{0\leq j\leq k}\binom{k}{j}W^{j-1}\lr{du-\f{x}{W}+\f{zd\bar d}{W}}^jx^{k-j}\al\\
&=\sls{0\leq j\leq k}\al_j\lr{du-\f{x}{W}+\f{zd\bar d}{W}}^j,
\end{align*}
where $\al_j=\binom{k}{j}W^{j-1}x^{k-j}\al$. Taking note that $M\leq D=X^\delta\ll\omega^C\frac{y}{W}$ for some constant $C$ whenever $\de<\f{\ta}{k}$, as well as $|\tau_3(d)\rd^2c_d|\leq\tau_3(d)$ for each $d\sim M$, we may apply \cite[Proposition 2.1]{MS} to obtain that there exists a positive integer $q\leq \omega^{-O_k(1)}$ such that for each $1\leq j\leq k$, $\|q\al_j\|\leq\omega^{-O_k(1)}\lr{\f{W}{y}}^j$, i.e.
\begin{align}\label{4.1}
\left\|\binom{k}{j}W^{j-1}x^{k-j}q\al\right\|\leq\omega^{-O_k(1)}\lr{\f{W}{y}}^j\ (1\leq j\leq k).
\end{align}
Now we suppose that $q'=[\binom{k}{1}q,\cdots,\binom{k}{k-1}W^{k-2}q,W^{k-1}q]$, then $q'\leq \omega^{-O_k(1)}W^k\leq\omega^{-O_k(1)}$, by recalling (\ref{w}) and $\omega\asymp(\log x)^{-A}$ again. Specially, when $j=k$ we can deduce from (\ref{4.1}) that
$$
\|q'\al\|\leq\|qW^{k-1}\al\|\leq\omega^{-O_k(1)}\lr{\f{W}{y}}^k.
$$
And for general $1\leq k-j\leq k$, if we assume that
\begin{align}\label{j}
	\|q'\al\|\leq\left\|\binom{k}{k-j}W^{k-j-1}q\al\right\|\leq\omega^{-O_k(1)}\f{W^k}{x^jy^{k-j}},
\end{align}
and consider the case of $k-(j+1)$. On the one hand, one can see from (\ref{4.1}) directly that
$$
\|q'x^{j+1}\al\|\leq\left\|\binom{k}{k-(j+1)}W^{k-j-2}x^{j+1}q\al\right\|\leq\omega^{-O_k(1)}\lr{\f{W}{y}}^{k-(j+1)}.
$$
On the other hand, we can deduce from the  hypothesis  (\ref{j}) that
$$
x^{j+1}\|q'\al\|\leq\omega^{-O_k(1)}\f{xW^k}{y^{k-j}}.
$$
On recalling that $y\asymp x^\ta$ and $\ta\in(1/2,1)$, above expression is less than $\f{1}{2}$ if and only if $k-j\geq2$, i.e. $j\leq k-2$. In total, on combining above discussion,  together with the fact that $x=\lfloor{X^{1/k}\rfloor}$ is an integer, one may find that when $j\leq k-2$,
$$
x^{j+1}\|q'\al\|=\|q'x^{j+1}\al\|\leq\omega^{-O_k(1)}\lr{\f{W}{y}}^{k-(j+1)}.
$$
In other word, when $j\leq k-2$ one could have
$$
\|q'\al\|\leq\omega^{-O_k(1)}\f{W^k}{x^{j+1}y^{k-(j+1)}}.
$$
In particular, when $j=k-2$ we have
$$
\|q'\al\|\leq\omega^{-O_k(1)}\f{W^k}{x^{k-1}y}.
$$
Recalling that  $\omega=(\log X)^{-A}$ and $W\ll\log\log x$ in (\ref{w}), the above conclusion can be written as---if (\ref{3}) holds, there exists $q\leq(\log X)^{c_kA}$ such that
$$
\|q\al\|\leq\f{(\log X)^{c_kA}}{x^{k-1}y}.
$$
Since in the case of $\al\in\mk$ either $q>(\log X)^{c_kA}$ or $\|q\al\|>\f{(\log X)^{c_kA}}{Y}$, this lemma therefore follows from $Y\asymp  X^{1-\f{1}{k}+\f{\ta}{k}}\asymp x^{k-1}y$.
\end{proof}
\end{Le}

\begin{Remark}
If the second result of \cite[Theorem 4.1]{Va} can be applied, similar minor arcs results should be obtained by making use of the methods in proving \cite[Lemma 4.2]{WW}. However, the appearance of W-trick makes the calculation rather complicated.	
\end{Remark}

On recalling (\ref{e}), we can rewrite $E_b(\al)$ as

\begin{align}\label{10}
E_b(\al)&=\sls{d_1,d_2|P}\rdo\rdt\sls{\f{X^{1/k}}{[d_1,d_2]}<t\leq \f{(X+Y)^{1/k}}{[d_1,d_2]}\atop [d_1,d_2]^kt^k\equiv b\md{W}}\ew{t^k[d_1,d_2]^k\al}\nonumber\\
&=\sls{d_1,d_2|P}\rdo\rdt f(b,[d_1,d_2],\al),
\end{align}
where $f(b,d,\alpha)$ is a function defined to be
\[
f(b,d,\alpha)=\sum_{\f{X^{1/k}}{d}<t\leq \f{(X+Y)^{1/k}}{d}\atop d^kt^k\equiv b\md{W}}e_W(d^kt^k\alpha).
\]

Now we define the exponential sum $V_q(a,b,d,c)$ as\begin{align}\label{vq}
V_q(a,b,d,c)=\sum_{z\in[W]\atop z^kd^k\equiv b\,(\textrm{mod}\,W)}\sum_{r(q)}e_{Wq}\big(a(z+Wr)^k+c(z+Wr)\big),	
\end{align}
then record \cite[Lemma 21]{Sal} as follows, since we'll use it frequently in the proof of major arcs.

\begin{Le}\label{lvq}
Let $a,b,d,q,k\in\mathbb N$ be such that $k\geq2$ and $(a,q)=(bd,W)=1$. Write $q=q_1q_2$ as $q_1$ is $w$-smooth and $(q_2,W)=1$. Then
\[
V_q(ad^k,b,d,0)=\xi(q)q_1\sum_{r(q_2)}e_{q_2}\big(a\psi_q(d)^k\bar{q_1W}r^k\big)\sum_{z(W)\atop z^k(d,q)^k\equiv b\,(\mathrm{mod}\,W)}\chi(z,(d,q)),
\]
where
\[
\psi_q(d)=\prod_{p^t\|d\atop p|q}p^t;
\]
\[
\bar{q_1W}q_1W\equiv 1\,(\mathrm{mod}\,q_2);
\]
\[
\chi(z,t)=e_{Wq}(at^kz^k)e_{q_2}(-at^kz^k\bar{q_1W});
\]
and
\[
\xi(q)=\begin{cases}
1\qquad&\text{if\ }q=1;\\
1&	\text{if\ }q_1|k\text{ and }q>w;\\
0&\text{otherwise.}
\end{cases}
\]	
\end{Le}

\begin{Le}[Major arc estimates]
Suppose that $1/2<\ta<1,k\geq2$ and $\al\in\mm$, then when $\de<\f{\ta}{k}$ we have
$$
|\widehat\nu_b(\al)-\widehat{1_{[N]}}(\al)|\ll N(\log\log\log N)^{-\frac{1}{2k}}.
$$
\begin{proof}
For frequency $\al\in\mm(q,a)$ with $1\leq a\leq q\leq Q$ and $(a,q)=1$, (specially, when $\alpha=0$ take $q=1$ and $a=0$) we may assume that $\al=\f{a}{q}+\beta$, thus, $|\beta|\leq\frac{Q}{qY}$. The application of \cite[Lemma 20]{Sal} gives that
\begin{align}\label{11}
f(b,[d_1,d_2],\al)=&\f{\vq}{qk[d_1,d_2]}X^{\f{1}{k}-1}\sls{X<t\leq X+Y\atop t\eq b\md{W}}\ew{\beta t}\nonumber\\
&+O\lr{(q,[d_1,d_2])W^2q^{\f{1}{2}+\va}\lr{1+\f{|\beta |Y}{W}}\log X+Y^2X^{\f{1}{k}-2}\f{1}{[d_1,d_2]}},	
\end{align}
where $V_q$ is defined in (\ref{vq}).

The next thing that needs to do is to provide a substitution of \cite[Lemma 22]{Sal} 
in the version of Selberg's upper bound sieve weights (\ref{r}), to be exact, to  bound the following expression
$$
\sls{d_1,d_2|P}\f{\rdo\rdt}{[d_1,d_2]}\vq.
$$

When $q=1$, then $a=0$ and $\beta=\al$, it is immediate from (\ref{vq}), the definition of $V_q$, that above representation becomes
\begin{align}\label{q1}
\sls{d_1,d_2|P}\f{\rdo\rdt}{[d_1,d_2]}\sls{z\in[W]\atop z^k[d_1,d_2]^k\eq b\md{W}}1=\sw\sls{d_1,d_2|P}\f{\rdo\rdt}{[d_1,d_2]},
\end{align}
as $\sw=\#\{z\in[W]:z^k\equiv b\pmod{W}\}$ and $[d_1,d_2]$ is coprime with $W$.

Now we assume that $q>1$ and write $q=q_1q_2$ where $q_1$ is $w$-smooth and $(q_2,W)=1$. When $q\neq1, q_1\nmid k$ or $q\leq w$, taking note that $\xi(q)=0$ in Lemma \ref{lvq}, we then deduce from  Lemma \ref{lvq} that
\begin{align}\label{q0}
\sls{d_1,d_2|P}\f{\rdo\rdt}{[d_1,d_2]}\vq=0.
\end{align}

It remains to handle the case when $q\neq1,q_1|k$ and $q>w$. Let us first classify the summation over $d_1$ and $d_2$ in the light of the greatest common divisor of $[d_1,d_2]$ and $q$,
$$
\sls{d_1,d_2|P}\f{\rdo\rdt}{[d_1,d_2]}\vq=\sls{t|q}\sls{\substack{
d_1,d_2|P\\
t|[d_1,d_2]\\
([d_1,d_2],\f{q}{t})=1 } }\f{\rdo\rdt}{[d_1,d_2]}\vq.
$$
Then another application of Lemma \ref{lvq} leads us to
\begin{align*}
&\qquad\sls{d_1,d_2|P}\f{\rdo\rdt}{[d_1,d_2]}\vq\\
&=\sls{t|q}\sls{\substack{
d_1,d_2|P\\
t|[d_1,d_2]\\
([d_1,d_2],\f{q}{t})=1 } }\f{\rdo\rdt}{[d_1,d_2]}q_1\sls{r \md{q_2}}e_{q_2}(a\psi_q([d_1,d_2])^k\overline{q_1W} r^k)\sls{z\in[W]\atop z^k([d_1,d_2],q)^k\eq b\md{W}}\chi(z,([d_1,d_2],q))\\
&=\sls{t|q}\sls{\substack{
d_1,d_2|P\\
t|[d_1,d_2]\\
([d_1,d_2],\f{q}{t})=1 } }\f{\rdo\rdt}{[d_1,d_2]}q_1\sls{r \md{q_2}}e_{q_2}(at^k\overline{q_1W} r^k)\sls{z \in[W]\atop z^kt^k\eq b\md{W}}\chi(z,t).
\end{align*}
Some comments for the last line may be required. The assumption $d_1,d_2|P$ yields that $d_1, d_2$ are square-free. And recalling $\psi_q(d)=\prod\limits_{p^t\|d\atop p|q}p^t=(d,q^\infty)$, we can deduce from $t=([d_1,d_2],q)$ and $[d_1,d_2]$ square-free that $\psi_q([d_1,d_2])=\psi_q(t)=t$. Since $(d_1d_2,W)=1, t|[d_1,d_2]$ and $q_1\leq k$, we see that
\begin{multline*}
	\qquad\sls{d_1,d_2|P}\f{\rdo\rdt}{[d_1,d_2]}\vq\\
\ll_k\sw\sls{t|q}\left|\sls{r\md{q_2}}e_{q_2}(at^k\overline{q_1W}r^k)\right|\cdot\left|\sls{d_1,d_2|P\atop t|[d_1,d_2]}\f{\rdo\rdt}{[d_1,d_2]}\right|.
\end{multline*}
By means of $(q_2,a\overline{q_1W})=1$, it is evident from \cite[Theorem 4.2]{Va} that
$$
\sls{r\md{q_2}}e_{q_2}(at^k\overline{q_1W}r^k)\ll q_2\lr{\f{(q_2,t^k)}{q_2}}^{\f{1}{k}}\ll
q_2^{1-\f{1}{k}}(q_2,t^k)^{\f{1}{k}}.
$$
We now recall that $\mathcal{J}=\sls{d|P\atop d<z}\f{1}{\phi(d)}=\sls{d<z\atop (d,W)=1}\f{|\mu(d)|}{\phi(d)}$, where $D=X^\delta$ is defined in (\ref{z}) and $P$ is the product of primes less than $D$. It can be deduced from \cite[Lemma A3]{Shao} that
$$
\sls{d_1,d_2|P\atop t|[d_1,d_2]}\f{\rdo\rdt}{[d_1,d_2]}\ll\mj t^{-1+\va}.
$$
Therefore, 
\begin{multline*}
\qquad\qquad\sls{d_1,d_2|P}\f{\rdo\rdt}{[d_1,d_2]}\vq\\
\ll \mj \sw q^{1-\f{1}{k}+\va}\sls{t|q}\f{(q,t^k)^{\f{1}{k}}}{t}
\ll\mj \sw q^{1-\f{1}{k}+\va},
\end{multline*}
as $\sls{t|q}\f{(q,t^k)^{\f{1}{k}}}{t}\leq\tau(q)\ll q^\va.$
 Besides, one can see from \cite[(A.1)]{Shao} that
\[
\mathcal J=\frac{\phi(W)}{W}(\log D+O_W(1)).
\]
We now combine above two expressions to conclude that that when $q\neq1,q_1|k$ and $q>w$,
\begin{align}\label{12}
\sls{d_1,d_2|P}\f{\rdo\rdt}{[d_1,d_2]}\vq\ll\f{\sw W}{\phi(W)}\f{q^{1-\f{1}{k}+\va}}{\log D}.	
\end{align}

Now let's turn our attention back to $\widehat\nu_b(\al)$. By taking note that 
$$
\sls{d_1,d_2|P}\rdo\rdt([d_1,d_2],q)\ll \sls{d\leq D}d^\va(d,q)\ll D^{1+\va},
$$
dues to $\rho_{d}$ is supported in the region $d\leq D$ and is also a 1-bounded function. Substituting (\ref{10}) and (\ref{11}) into (\ref{v}) yields that
\begin{align*}
\widehat\nu_b(\al)&=e\lr{-(\f{b}{W}+m)\al}\f{\phi(W)}{\al^+W \sw} \f{\log X}{qk}\sls{X<t\leq X+Y\atop t\eq b\md{W}}\ew{\beta t}\sls{d_1,d_2|P}\f{\rdo\rdt}{[d_1,d_2]}\\
&\cdot\vq+O\lr{D^{1+\va}X^{1-\f{1}{k}}W^2q^{\f{1}{2}+\va}\lr{1+\f{|\beta|Y}{W}}(\log X)^2+Y^2X^{-1}(\log X)^2}.
\end{align*}
Let us deal firstly the above $O$-term. To begin with, the facts that(\ref{yy}), (\ref{xy})  and $\ta$ is strictly less than $1$ lead us to
$$
Y^2X^{-1}(\log X)^2\ll X^{1-\f{2(1-\theta)}{k}}(\log X)^2\ll\frac{N}{\log N}.
$$
Secondly, one deduces from $q\leq Q$, $|\beta|\leq\f{Q}{qY}$, together with (\ref{w}), (\ref{z}) and (\ref{q}) that
\begin{multline*}
D^{1+\va}X^{1-\f{1}{k}}W^2q^{\f{1}{2}+\va}\lr{1+\f{|\beta|Y}{W}}(\log X)^2
\ll X^{1-\f{1}{k}+\de+\va}W(WQ^{1/2+\va}+Q^{1+\va})\\
\ll NX^{\de-\f{\ta}{k}+\va}W(WQ^{1/2+\va}+Q^{1+\va}),
\end{multline*}
thus when $\de<\f{\ta}{k}$ and $0<\va<1$ is arbitrarily small, it is bounded by $O(\frac{N}{\log N})$. We then have
\begin{align}\label{mian}
\widehat\nu_b(\al)&=e\lr{-(\f{b}{W}+m)\al}\f{\phi(W)}{\al^+W \sw} \f{\log X}{qk}\sls{X<t\leq X+Y\atop t\eq b\md{W}}\ew{\beta t}\nonumber\\
&\cdot\sls{d_1,d_2|P}\f{\rdo\rdt}{[d_1,d_2]}\vq+O(\frac{N}{\log N}).	
\end{align}

Suppose that $q>1$, if $q_1\nmid k$ or $q\leq w$, substituting  (\ref{q0}) into (\ref{mian}) will lead to $|\widehat\nu_b(\alpha)|\ll\frac{N}{\log N}$; if $q_1|k$ and $q>w$, substituting  (\ref{12}) into (\ref{mian}) gives that
\begin{multline*}
\widehat\nu_b(\al)\ll_{\va,k}\bigg|\sls{X<t\leq X+Y\atop t\eq b\md{W}}\ew{\beta t}\bigg|q^{\va-\f{1}{k}}+O(\frac{N}{\log N})\\
\ll w^{\epsilon-\frac{1}{k}}N+O(\frac{N}{\log N})\ll N(\log\log\log N)^{-\frac{1}{2k}}.
\end{multline*}
When $q=1$, for $\alpha\in\mathcal M(q)$, we must have $a=0$ and $\alpha=\beta\approx 0$, so
\[
\widehat\nu_b(\alpha)=\widehat\nu_b(0)+O(N/\log N).
\] 
Whilst, by substituting (\ref{q1}) into (\ref{mian}) one has
\begin{align*}
\widehat\nu_b(0)&=\f{\phi(W)}{\al^+W \sw} \f{\log X}{k}\sls{d_1,d_2|P}\f{\rdo\rdt}{[d_1,d_2]}V_q(0,b,[d_1,d_2],0)\sls{X<t\leq X+Y\atop t\eq b\md{W}}1+O(\frac{N}{\log N})\\
&=\f{\phi(W)}{\al^+W } \f{\log X}{k}\sls{d_1,d_2|P}\f{\rdo\rdt}{[d_1,d_2]}\sls{X<t\leq X+Y\atop t\eq b\md{W}}1+O(\frac{N}{\log N}).
\end{align*}
Recalling that $\al^+=\f{\phi(W)}{kW}\log X\sls{d_1,d_2|P}\f{\rdo\rdt}{[d_1,d_2]}$ which is defined in (\ref{1}) as well as $Y=WN$ in (\ref{xy}), on combining the above two expressions, one has
\[
\widehat\nu_b(\alpha)=\sls{X<t\leq X+Y\atop t\eq b\md{W}}1+O(\frac{N}{\log N})=N+O(\frac{N}{\log N}).
\]

Next, we can deal with $\widehat{1_{[N]}}(\al)$ with $\alpha\in\mathcal M$. Noting that $\alpha\neq0$ whenever $q>1$, and thus,
\[
\|\alpha\|\geq\frac{1}{q}-\big|\alpha-\frac{a}{q}\big|\geq \frac{1}{q}-\frac{Q}{qY}\gg \frac{1}{q}.
\]
Hence, when $\alpha\in\mathcal M$ and $q>1$, it follows from Fourier transform that
\[
|\widehat{1_{[N]}}(\al)|=\Big|\frac{1-e((N+1)\alpha)}{1-e(\alpha)}\Big|\ll\|\alpha\|^{-1}\ll\frac{N}{\log N}.
\]
When $q=1$ and $\alpha=0$ we plainly have $\widehat{1_{[N]}}(0)=N$. To sum up, when $\alpha\in\mathcal M$,  we always have
$$
\widehat\nu_b(\al)=\widehat{1_{[N]}}(\al)+O(N(\log\log\log N)^{-\frac{1}{2k}}).
$$

\end{proof}

\end{Le}

By combining Lemma 5 and Lemma 7 one finds that when $k\geq2$,   $\ta\in(\f{1}{2},1)$ and $\de<\f{\ta}{k}$, we would have
$$
|\widehat\nu_b(\al)-\widehat{1_{[N]}}(\al)|\ll N(\log\log\log N)^{-\frac{1}{2k}},
$$
for every $\al\in\mt$. And this is Lemma 4.

\vspace{2mm}

\section{An almost all version of transference principle}

We would like to establish an almost all version of transference principle, and then apply it to our problem in next section. The below result is a generalization of \cite[Lemma 5]{Sal}, and the proof is just a refinement of Salmensuu's argument, so we record it in the appendix.

\begin{Le}\label{1bou}
Suppose that $s\geq3$ is an integer, $0<\tilde\kappa<1$ and $\frac{1}{s}\leq\alpha<1-\tilde\kappa$. Then there exists a number $\eta=\eta(\tilde\kappa,s)\in(0,1)$ such that the following property holds.

Let $g_1,\dots,g_s:[N]\to[0,1+\tilde\kappa/4]$ be non-negative functions such that $\E_P(g_i)\geq\alpha+\tilde\kappa$ for every arithmetic progression $P\subseteq[N]$ with $|P|\geq\eta N$. Then for all $n\in[N/2,N]$, we have
\[
g_1*\dots*g_s(n)\gg_{s,\tilde\kappa}N^{s-1}.
\]
\end{Le}

The next result is the so-called almost all version transference principle. The argument is that for any unbounded function, we can always decompose it into an essentially bounded component and a small component (in the sense of its Fourier coefficients are small). Then the dominants are from the bounded functions, and this will be proven by using  Lemma \ref{1bou}. Whilst one can also show that, on average, the contribution contains at least one small component is negligible.

\begin{Proposition}\label{almost}
Let $N$ be a sufficiently large number. Let $f_1,\cdots,f_s:[N]\rightarrow\mathbb{R}_{\geq0}$ be non-negative functions. And each of $f\in\{f_1,\cdots,f_s\}$ satisfies the following three conditions:
  \subsection*{1)(Mean Condition)} There exist numbers $\tilde\kappa>0$ and $\eta\in(0,1)$ such that 
   $$\me_{n\in P}f(n)\geq\f{1}{s}+\tilde\kappa$$
   holds for every arithmetic progression $P\subset[N]$ with $|P|\geq\eta N$;
 \subsection*{2)(Pseudorandomness Condition)}There is a function $\Delta(N)$ which tends to zero as $N$ tends to infinity, and a majorant function $\nu:[N]\rightarrow\mathbb{R}_{\geq0}$ such that $f\leq\nu$ pointwise and 
 \[
 \|\widehat{\nu}-\widehat{1_{[N]}}\|_{\infty}\leq\Delta(N) N;
 \]	 
\subsection*{3)(Restriction Estimate)} There exists  $q\in(2s-1,2s)$ such that $\|\widehat{f}\|_q^q\ll N^{q-1}$.

Then the inequality
$$
f_1*\cdots*f_s(n)\gg_{s,\tilde\kappa,q} N^{s-1}
$$	
fails for at most $O\Bigbrac{\frac{N}{|\log \Delta(N)|^{\frac{2s-q}{6s+3}}}}$ values of $n\in[N/2,N]$.

\begin{proof}

Let $0<\va_1<1$ be a parameter to be specify later. Define the large spectrum set of $f$ as
\[
\Omega=\bigset{\alpha\in\T:|\hat f(\alpha)|\geq\va_1 N},
\]
and the Bohr set as
\[
B=\bigset{|n|\leq\va_1 N:\norm{n\alpha}\leq\va_1\text{ for all }\alpha\in\Omega}.
\]
\cite[Lemma A.2]{Pre} gives that
\begin{align}\label{bohr}
|B|\geq\va_1^{O_q(\va_1^{-q-1})}N=:\epsilon_2N	,
\end{align}
where $q\in(2s-1,2s)$ is the parameter in condition 3. We now define a ``smooth" version of $f$ together with  their difference as follows,
\[
g(n)=\E_{x_1,x_2\in B}f(n+x_1-x_2)\qquad h=f-g.
\]
For convenience, write $\sigma=|B|^{-1}1_B$ as the normalized indicator function of $B$, and let $\sigma_1(n)=\sigma(-n)$. Then $g=f*\sigma*\sigma_1$. Besides, one can also find from Fourier inverse formula, as well as $f\leq\nu$ pointwise, that
\[
	g(n)\leq\int_\T\hat{\nu*\sigma*\sigma_1}(\gamma)e(-n\gamma)\mathrm d\gamma=\int_\T\hat\nu(\gamma)|\hat\sigma(\gamma)|^2e(-n\gamma)\mathrm d\gamma.
\]
On the other hand, it can be deduced from condition 2, Parseval's identity, together with (\ref{bohr}) that
	\begin{align*}
	\int_\T|\hat\nu(\gamma)-\hat1_{[N]}(\gamma)||\hat\sigma(\gamma)|^2e(-n\gamma)\mathrm d\gamma &\leq\|\hat\nu-\hat1_{[N]}\|_\infty\int_\T|\hat\sigma(\gamma)|^2\mathrm d\gamma\\
	&\leq\Delta( N)N\norm{\sigma}_2^2\leq\Delta( N)N|B|^{-1}
	\leq\Delta( N)\epsilon_2^{-1}.
	\end{align*}
On combining above two inequalities and then applying triangle inequality, one has
	\begin{align}\label{g1}
	g(n)&\leq\int_\T\hat1_{[N]}(\gamma)|\hat\sigma(\gamma)|^2e(-n\gamma)\mathrm d\gamma+\Delta( N)\epsilon_2^{-1}\nonumber\\
	&=\E_{b_1,b_2\in B}1_{[N]}(n+b_1-b_2)+\Delta( N)\epsilon_2^{-1}\leq 1+\Delta( N)\epsilon_2^{-1}.	\end{align} 
	
	Let $\eta(\frac{2}{3}\tilde\kappa,s)$ be the quantity such that Lemma \ref{1bou} holds for parameters $2\tilde\kappa/3,s$, and then take $\eta=\eta(\frac{2}{3}\tilde\kappa,s)/2$ in condition 1. Now we suppose that $P\subset[N]$ is an arithmetic progression with $|P|\geq\eta(\frac{2}{3}\tilde\kappa,s) N=2\eta N$, from condition 1, one has,
	\[
	\E_{n\in P}g(n)=\E_{b_1,b_2\in B}|P|^{-1}\sum_{n\in (P+b_1-b_2)\cap[N]}f(n)\geq\frac{1}{s}+\tilde\kappa-\epsilon_1/\eta\geq\frac{1}{s}+\frac{2\tilde\kappa}{3},
	\]
	whenever $\epsilon_1/\eta<\frac{\tilde\kappa}{3}$. Recalling that $2s-1<q<2s$, when we take $\epsilon_1=|\log \Delta(N)|^{-\frac{1}{2s+1}}$, we would have both $\Delta(N)\epsilon_2\leq\tilde\kappa/10$ and $\epsilon_1/\eta<\frac{\tilde\kappa}{3}$. Thus, an application of Lemma \ref{1bou} to $g_1,\cdots,g_s$ yields that for every $n\sim N/2$,
 $$
 g_1*\cdots*g_s(n)\gg_{\tilde\kappa,s}N^{s-1}.
 $$

It therefore suffices to show that, there is a parameter $\epsilon'=\epsilon'(\epsilon_1,s,q)$ to be chosen later, such that out of at most $\epsilon'N$ of $n\sim N/2$, we have
\begin{align}\label{claim}
y_1*\cdots*y_s(n)\ll\va' N^{s-1},
\end{align}
where $y_i$ is either $g_i$ or $h_i$, and at least one $y_i$ is equal to $h_i$. 

We need to do some preparation before we prove  (\ref{claim}). For every frequency $\gamma\in\T$,
 	\[
	\hat h(\gamma)=\hat f(\gamma)-\hat g(\gamma)=\hat f(\gamma)(1-|\hat\sigma(\gamma)|^2).
	\]
When $\gamma\not\in\Omega$, $|\hat h(\gamma)|\leq{\va_1} N$ as $0\leq 1-|\hat\sigma(\gamma)|^2\leq1$. When $\gamma\in\Omega$, it follows from triangle inequality, condition 2, as well as $\|\hat f\|_\infty=\sum_n f(n)\leq\sum_n\nu(n)$ that
	\[
	\|\hat f\|_\infty\leq\sum_{n}1_{[N]}(n)+\Bigabs{\sum_{n}(1_{[N]}(n)-\nu(n))}\leq(1+\Delta(N))N.
	\]
	Thus, when $\gamma\in\Omega$ we have
\[
|\hat h(\gamma)|\leq(1+\Delta(N))N(1-|\hat\sigma(\gamma)|^2).
\]
Whilst, from the definition of Bohr set $B$, we also have
\[
\bigabs{1-|\hat\sigma(\gamma)|^2}\ll\E_{b\in B}\norm{b\gamma}\ll{\va_1}.
\] 
Taking note that $\Delta(N)\to0$ as $N\to\infty$, whether $\gamma\in\Omega$ or not, we always have
\begin{align}\label{h}
\|\widehat{h}\|_{\infty}	\ll{\va_1} N.
\end{align}

Now we back our attention to prove (\ref{claim}).  Without loss of generality, we can assume that $y_1=h_1$. Let $\mathcal E$ be the set of $n\sim N/2$ for which if $n\in\mathcal E$ then
\[
|h_1*\cdots*y_s(n)|\gg\va'N^{s-1}.
\] 
On summing over $n\in\mathcal E$ and using Cauchy-Schwarz inequality to get that
\begin{align*}
\va'|\mathcal E|N^{s-1}\ll
\sls{n\in\mathcal E}|h_1*\cdots*y_s(n)|&\leq|\mathcal E|^{\f{1}{2}}\left(\sls{n\in\mathcal E}\big|h_1*\cdots*y_s(n)\big|^2\right)^{\f{1}{2}}\\
&\leq|\mathcal E|^{\f{1}{2}}\left(\sls{n\in[N]}\big|h_1*\cdots*y_s(n)\big|^2\right)^{\f{1}{2}}.
\end{align*}
The making use of Plancherel formula leads us to
\begin{align*}
\sls{n\in[N]}\big|h_1*\cdots*y_s(n)\big|^2&=\int_{\mt}|\widehat{h_1*\cdots*y_s}(\al)|^2d\al\\
&=\int_\mt|\widehat{h_1}(\al)\cdots\widehat{y_s}(\al)|^2d\al\\
&\leq \|\widehat{h_1}\|_{\infty}^{2s-q}\int_\mt|\widehat{h_1}(\al)|^{2-(2s-q)}|\widehat{y_2}(\al)|^2\cdots|\widehat{y_s}(\al)|^2d\al.
\end{align*}
By combining the above two expressions and then making use of H\"{o}lder's inequality together with  (\ref{h}),  the restriction estimate (condition 3), also the fact that $|\hat g(\gamma)|,|\hat h(\gamma)|\leq|\hat f(\gamma)|$ pointwise, we would have
\begin{align*}
(\va'|\mathcal E|N^{s-1})^2&\leq |\mathcal E|\|\widehat{h_1}\|_{\infty}^{2s-q}\int_\mt|\widehat{h_1}(\al)|^{2-(2s-q)}|\widehat{y_2}(\al)|^2\cdots|\widehat{y_s}(\al)|^2d\al\\
&\leq|\mathcal E|\|\widehat{h_1}\|_{\infty}^{2s-q}
\max_{i}\|\widehat{f_i}\|_q^{q}\ll \va_1^{2s-q}N^{2s-1}|\mathcal E|.
\end{align*}
Thus, $|\mathcal E|\ll\va_1^{\frac{2s-q}{3}}N$  by setting $\epsilon'=\epsilon_1^{\frac{2s-q}{3}}$. And the proposition follows from $\epsilon_1=|\log \Delta(N)|^{-\frac{1}{2s+1}}$.

\end{proof}

\end{Proposition}

\bigskip

\section{Proof of Theorem 3}

\begin{Le}\label{residue}
Let $k\geq2$, $s>\f{k(k+1)}{2}$ and $h$ be positive integers. Let $R_k$ be the integer defined in (\ref{rk}) and $(h,R_k)$ be the greatest common divisor of $h$ and $R_k$.  Let $m\in\mathbb{Z}_h$ such that $m\eq s\pmod{(h,R_k)}$, particularly, when $k=3$ and $s=7$, $m$ satisfies an additional condition $m\not\equiv0\pmod{9}$. Define
\begin{align*}
M_s(h,m)=\#\{(y_1,\cdots,y_s)\in\mathbb{Z}_h^*\times\cdots\times\mathbb{Z}_h^*:y_1^k+\cdots+y_s^k\eq m\pmod{h}\}.
\end{align*}
Then we have $M_s(h,m)>0$.

\begin{proof}
It can be shown that $M_s(h,m)$ is multiplicative as a function of $h$. So we can just consider the case of $h=p^t$. According to the proof of \cite[Lemma 10]{Sal} when $t\geq\gamma=\gamma(k,p)$, where $\gamma$ is defined in (\ref{ga}), one could have
$$
p^tM_s(p^t,m)=p^{s(t-\gamma)}M_s(p^\gamma,m).
$$
We would continue the proof by distinguishing into two cases according to whether $p-1$ and $k$ are coprime or not.

\subsection*{(case 1).}When $p-1|k$, one can deduce from (\ref{rk}) that $(p^\gamma,R_k)=p^\gamma$, and then the congruent condition $m\eq s\pmod{(p^\gamma,R_k)}$ is equal to  $m\eq s\pmod{p^\gamma}$ . Thus in this case $M_s(p^\gamma,m)>0$ follows from \cite[Lemma 8.9]{Hua}.
\subsection*{(case 2).}  The case of $p-1\nmid k$ will be discussed on several subcases according to the value of $k$.
 \subsection*{(cses 2.1).} When $k\geq5$ and $s>\f{k(k+1)}{2}$, we always have $s\geq3k$. Hence, $M_s(p^\gamma,m)>0$ follows from the first result of \cite[Lemma 8.8]{Hua};
\subsection*{(case 2.2).} When $k=4$, for any pirme $p$ with $p-1\nmid k$, one finds that $8\neq p^\tau(p-1)$ always holds, where $\tau$ is the number such that $p^\tau\|k$. And $M_s(p^\gamma,m)>0$ may be established by the second result of \cite[Lemma 8.8]{Hua} together with the fact that $s>\f{k(k+1)}{2}\geq2 k$;
\subsection*{(case 2.3)} When $k=3$ and $p\neq3,7$, similar to case 2.2 we can get $M_s(p^\gamma,m)>0$. If $p=3$ or $p=7$, we can still obtain $M_s(p^\gamma,m)>0$ by making use of the third result of \cite[Lemma 8.8]{Hua};
\subsection*{(case 2.4)} When $k=2$, it suffices to show that for any prime $p\neq2,3$ and any $m\in\mz_p$ the following quadratic congruent equation has a solution of $y_1,\cdots,y_4\in\mz_p^*$
$$
y_1^2+\cdots+y_4^2\eq m\pmod{p}.
$$
Define a set $\my$ by setting $\my=\{y^2:y\in\mz_p^*\}$. Then $|\my|=\f{p-1}{2}$. And applying a result of Cauchy-Davenport-Chowla (see \cite[Lemma 2.14]{Va} or \cite[Lemma 8.7]{Hua} as an example) one concludes that 
$$
|\my+\my+\my+\my|\geq p,
$$
 whenever prime $p>3$. Thus, the above quadratic congruent equation is solvable for all $m\in\mz_p$.

The case of $t<\gamma(k,p)$ will only occur in $p=2$, and then $(R_k,2^t)=2^t$ takes us from the congruent condition $m\eq s\pmod{(R_k,2^t)}$ to $m\eq s\pmod{2^t}$. Obviously, $y_1=\cdots=y_s\eq1\pmod{2^t}$ is a solution of $y_1^2+\cdots+y_s^2\eq m\pmod{2^t}$. Combine all of above situations we complete the proof this lemma.
\end{proof}

\end{Le}

\bigskip

We are going to continue using the assumptions of $s,k,W,m,N$ that appeared in Section 3 at the moment. And suppose that $1\leq b_1,\dots, b_s<W$ are integers with $(b_1\dots b_s,W)=1$. Using these integers, we can define a group of functions $f_{b_i}:[N]\to\R_{\geq0}$ in the following way 
\begin{align}\label{fi}
f_{b_i}(n)=	
\begin{cases}
\f{\phi(W)}{\al^{+}W\sigma_W(b_i)}X^{1-\f{1}{k}}\log X,\qquad &\text{if\ }W(n+m)+b_i=p^k;\\
0,&\text{otherwise,}
\end{cases}
\end{align}
and also define the corresponding majorants by taking $\nu_b=\nu_{b_i}$ in (\ref{vb}).

We then hope the application of Proposition \ref{almost} would lead us to the desired result. However, it seems hard to do so. Because Proposition \ref{almost} is not a pointwise result, instead, it just gives us the information of an almost all type result in a given interval. To make Proposition \ref{almost} work, a possible way is to divide the interval $[1,M]$ into subsets for which elements $(\,\textrm{mod}\, W\,)$ are the same. Then the making use of Proposition \ref{almost} would conclude that the density of the subset of $\set{n\leq M:n\equiv b\mod{W},}$ (with $(b,W)=1$) for which elements cannot be a sum of $k$-th power of primes is $O(|\log\Delta(N)|^{-c})$. However, dues to our $\Delta(N)$ in Lemma 4 is too large compared with $\phi(W)^{-1}$, the probability that a number in a specifically reduced residue class modulus $W$, this fails us to get the desired result. To overcome this, we'd like to modify Proposition \ref{almost} such that it is not so depending on the congruent condition, and then wish this argument would be sufficient to prove our Theorem 3.

\vspace{3mm}

\noindent\emph{Proof of Theorem 3.}

\vspace{2mm}

Suppose that
\begin{align}\label{comf}
f=\sum_{1\leq b\leq W\atop(b,W)=1}f_b=:\mathop{{\sum}^*}_{b(W)}f_b,	
\end{align}
then, from (\ref{fi}), $f(n)$ is non-zero if there is a prime $p$ such that $W(m+n)+b=p^k$ for some integer $(b,W)=1$.

\subsection*{Cliam} If for arbitrary $1\leq b_1,\dots,b_s\leq W$ with $(b_1\dots b_s,W)=1$, $f_{b_1},\dots,f_{b_s}$ satisfies the three conditions in Proposition \ref{almost}, then
\begin{align}\label{fcon}
\begin{matrix}
\underbrace{f*\cdots*f}(n) &\gg_{s,\tilde\kappa,q} \phi(W)^sN^{s-1} \\ 
 s \text{-fold }
& 
\end{matrix}
\end{align}	
fails for at most $O\Bigbrac{\frac{N}{|\log \Delta(N)|^{\frac{2s-q}{6s+3}}}}$ values of $n\in[N/2,N]$.

We just have a quick look at how to prove this claim here. Firstly as one may see from Proposition \ref{almost} that 
\[
g_{b_1}*\cdots*g_{b_s}(n)\gg_{\tilde\kappa,s}N^{s-1}
\]
holds for every $n\sim N/2$. Thus, for every $n\sim N/2$, we have
\[
\Bigbrac{\mathop{{\sum}^*}_{b_1(W)}g_{b_1}}*\cdots*\Bigbrac{\mathop{{\sum}^*}_{b_s(W)}g_{b_s}}(n)\gg\phi(W)^s\min_{b_1,\dots,b_s}g_{b_1}*\cdots*g_{b_s}(n)\gg_{\tilde\kappa,s}\phi(W)^sN^{s-1}.
\]
Secondly, define the set
\[
\mathcal E=\Bigset{n\sim N/2:\Bigabs{\Bigbrac{\mathop{{\sum}^*}_{b_1(W)}h_{b_1}}*\Bigbrac{\mathop{{\sum}^*}_{b_2(W)}f_{b_2}}*\cdots*\Bigbrac{\mathop{{\sum}^*}_{b_s(W)}f_{b_s}}(n)}\gg_{s,\tilde\kappa}|\log\Delta(N)|^{-\frac{2s-q}{6s+3}}\phi(W)^sN^{s-1}}.
\]
We would have
\begin{multline*}
|\log\Delta(N)|^{-\frac{2s-q}{6s+3}}\phi(W)^sN^{s-1}	|\mathcal E|\ll|\mathcal E|^{1/2}\biggbrac{\sum_{n\in[N]}\Bigabs{\Bigbrac{\mathop{{\sum}^*}_{b_1(W)}h_{b_1}}*\cdots*\Bigbrac{\mathop{{\sum}^*}_{b_s(W)}f_{b_s}}(n)}^2}^{1/2}\\
\leq|\mathcal E|^{1/2}\Bigbrac{\phi(W)^{2s}\sum_{n\in[N]}\max_{b_1,\dots,b_s}\bigabs{h_{b_1}*\cdots*f_{b_s}(n)}^2}^{1/2}.
\end{multline*}
That is,
\[
|\log\Delta(N)|^{-\frac{4s-2q}{6s+3}}|\mathcal E|N^{2s-2}\ll\sum_{n\in[N]}\sup_{b_1,\dots,b_s}\bigabs{h_{b_1}*\cdots*f_{b_s}(n)}^2.
\]
Thus, one can find this claim follows from the assumption that Proposition \ref{almost} is valid for any choice of $(b_1\dots b_s,W)=1$.

We now partition the interval $[1,N]$ into $[1,N|\log\Delta(N)|^{-\frac{2s-q}{6s+3}}]$ and dyadic intervals $(\frac{N}{2^i},\frac{N}{2^{i-1}}]$ with $1\leq i\ll|\log\log \Delta(N)|$. For a fixed number $1\leq i\ll|\log\log \Delta(N)|$, one can see from the claim that except at most $O\bigbrac{{|\log \Delta(N)|^{-\frac{2s-q}{6s+3}}}\frac{N}{2^{i-1}}}$ of $n\sim\frac{N}{2^i}$,
\[
\begin{matrix}
\underbrace{f*\cdots*f}(n) &\gg_{s,\tilde\kappa,q} \phi(W)^s\bigbrac{\frac{N}{2^{i-1}}}^{s-1} \\ 
 s \text{-fold }
& 
\end{matrix}.
\]	
The summing over the intervals yields that the inequality (\ref{fcon}) fails for at most 
\[
N|\log\Delta(N)|^{-\frac{2s-q}{6s+3}}+\sum_{1\leq i\ll|\log\log \Delta(N)|}{|\log \Delta(N)|^{-\frac{2s-q}{6s+3}}}\frac{N}{2^{i-1}}\ll N|\log\Delta(N)|^{-\frac{2s-q}{6s+3}}
\]
of $n\in[N]$.

We are planning to verify the three conditions in Proposition \ref{almost} for the functions $f_{b_1},\dots,f_{b_s}$. The pseudorandomness condition in Proposition \ref{almost} can be verified by Lemma 4 and by taking $\Delta(N)=(\log\log\log N)^{-1/2k}$. \cite[Lemma 7]{Sal}  also shows that when $\f{\al^-}{\al^+}(1-\va)>\f{1}{s}$, we have $\me_{n\in P}f_{b_i}(n)\geq\f{1}{s}+\va$ for every long arithmetic progression $P$. Dues to the definition of $\alpha^+$ in (\ref{al}) we could have $\f{\al^-}{\al^+}>\f{\al^-k\delta}{2}$, and as a consequence, the mean value condition holds whenever $s>\f{2}{\al^-k\delta}$. Moreover, \cite[Lemma 9]{Sal} shows that the restriction lemma holds for every number $q\in(2s-1,2s)$ whenever $2s>k(k+1)$. To sum up all above analysis and setting $q=2s-1/2$, we have the following conclusion.

\begin{Corollary} \label{coro}
	Let $\al^->0$ be the parameter in Theorem 1. Suppose that $f:[N]\to\C$ is the function defined in (\ref{comf}). Then when $s>\max\{\f{2}{\al^-\ta},\f{k(k+1)}{2}\}$ the following inequality 
		\[
\begin{matrix}
\underbrace{f*\cdots*f}(n) &\gg_{s} \phi(W)^{s}N^{s-1} \\ 
 s \text{-fold }
& 
\end{matrix}
\]
fails for at most $O\bigbrac{(\log_4 N)^{-\frac{1}{12s+6}}N}$ of $n\in[N]$.
\end{Corollary}

We are now ready to prove Theorem 3. Indeed, Theorem 3 may be derived from the following statement. Let $M$ be a sufficiently large integer. All but at most $O\bigbrac{(\log_4 M)^{-\frac{1}{12s+6}} M}$ of integers $n\equiv s\pmod{R_k}$ lie in the interval $(\f{M}{2},M]$, there are primes $p_1,\cdots,p_s\in\left((\f{n}{s})^{\f{1}{k}}-n^{\f{\ta}{k}},(\f{n}{s})^{\f{1}{k}}+n^{\f{\ta}{k}}\right]$ such that $n=p_1^k+\cdots+p_s^k$. To begin with, we divide the interval $(\f{M}{2},M]$ into subintervals $(M_i-(\f{M_i}{s})^{\f{k-1+\ta}{k}},M_i](i\geq1)$ with $M_1=M$ and $M_{i+1}=M_i-(\f{M_i}{s})^{\f{k-1+\ta}{k}}$. Then the number of such subintervals is $i\ll M^{\frac{1-\theta}{k}}$. Since the treatment for each of these intervals is the same, we now focus on one such interval, and write for abbreviation as $(M-H,M]$ with $H=\bigbrac{\frac{M}{s}}^{\frac{k-1+\theta}{k}}$.

Let $x=(\f{M}{s})^{\f{1}{k}}$, and take 
$$
N=\Bigfloor{\f{(x-x^\ta-W)^k-(x-\f{x^\ta}{s})^k}{W}}\qquad\text{and}\qquad m=\Bigfloor{\f{(x-\f{x^\ta}{sk})^k}{W}}.
$$
We see that
\begin{align}
\label{21}
N&=W^{-1}(1+o(1))\f{ks+k}{s}x^{k-1+\ta}\\
Y&=WN=(1+o(1))\f{ks+k}{s}x^{k-1+\ta},	\nonumber\\
\label{22}
&Wm=x^k-(1+o(1))\f{1}{s}x^{k-1+\ta}.
\end{align}
Corollary \ref{coro}, as well as the definition of convolution, tells us that out of at most $O\bigbrac{(\log_4 N)^{-\frac{1}{12k+6}}N}$ of $n\in[N]$,
\[
\sum_{n_1+\dots+n_s=n}f(n_1)\dots f(n_s)>0.
\]
From the definition of the function $f$ in (\ref{comf}), there are primes $p_1,\dots,p_s$ and integers $(b_1\dots b_s,W)=1$ such that
\[
sWm+Wn_0+b_1+\dots+b_s=p_1^k+\dots+p_s^k.
\]
Besides, it can be deduced from (\ref{21}) and (\ref{22}), together with $1\leq b_1,\dots,b_s\leq W$ that sequence $\set{sWm+Wn_0+b_1+\dots+b_s}_{n_0\leq N,(b_1\dots b_s,W)=1}$ lies in the interval $(M-H,M+O(W)]$. Therefore, in view of Lemma \ref{residue}, there are at most $O\bigbrac{(\log_4 M)^{-\frac{1}{12s+6}}H}$ of $n\in(M-H,M]$ with $n\equiv s\pmod{R_k}$ such that the following equation fails
\[
n=p_1^k+\dots+p_s^k.
\]
We then complete the proof by direct checking that $|p_i^k-\frac{M}{s}|\ll\bigbrac{\frac{M}{s}}^{1-\frac{1}{k}+\frac{\theta}{k}}$.

\qed

\appendix

\section{Proof of Lemma \ref{1bou}}

We are going to prove Lemma \ref{1bou} in this appendix. As we said before the proof is similar to \cite[Lemma 5]{Sal}, so we hope this would serve as an easy-read exposition.

Let $A$ and $B$ be two subsets of $[N]$. $0<\eta<1$ is a parameter, we define the  \textit{$\eta$-popular sumset} of sets $A$ and $B$ as follows,
\[
S_\eta(A,B)=\set{n:1_A*1_B(n)\geq\eta\max\set{|A|,|B|}}.
\]
\begin{Le}\label{sum}
Suppose that $\tilde\kappa,\delta\in(0,1)$, there is a number $\eta=\eta(\tilde\kappa,\delta)>0$ such that the following statement holds. Let $\alpha,\beta\in (\tilde\kappa,1)$, and $A,B\subseteq[N]$ satisfy $|A\cap P|\geq\alpha|P|$ and $|B\cap P|\geq\beta|P|$	for all progressions $P\subseteq[N]$ with $|P|\geq\eta N$. Then for any progression $Q\subseteq[2N]$ of length $|Q|\geq2\delta N$, we have
\[
|S_\eta(A,B)\cap Q|\geq\Bigbrac{\min\set{\alpha+\beta,1}-\tilde\kappa}|Q|.
\]
\begin{proof}
See \cite[Lemma 4]{Sal}.	
\end{proof}
	
\end{Le}

We would like to apply above lemma repeatedly to prove Lemma \ref{1bou}.

\vspace{3mm}

\noindent\emph{Proof of Lemma \ref{1bou}.}

\vspace{2mm}
Fixed a number $n_0\in[N/2,N]$, set $N_1=\floor{\frac{n_0}{2^{s-2}}}$,\ $N_{i+1}=2N_i(1\leq i\leq s-2)$, and define
\[
A_1=\set{n\in[N_1]:g_1(n)\geq\tilde\kappa/4},\quad A_i=\set{n\in[N_{i-1}]:g_i(n)\geq\tilde\kappa/4}(2\leq i\leq s).
\]
Since $g_i$ is non-negative, it is clear that $g_i(n)\geq\frac{\tilde\kappa}{4}1_{A_i}(n)$, and thus,
\[
g_1*\dots*g_s(n_0)\geq\bigbrac{\frac{\tilde\kappa}{4}1_{A_1}}*\dots*\bigbrac{\frac{\tilde\kappa}{4}1_{A_s}}(n_0)\gg\tilde\kappa^s1_{A_1}*\dots*1_{A_s}(n_0).
\]
It then suffices to prove that  for above $n_0\in[N/2,N]$ we have
\[
1_{A_1}*\dots*1_{A_s}(n_0)\gg_{s,\tilde\kappa}N^{s-1}.
\]
It is convenient to take $\kappa'=\frac{\tilde\kappa}{4s}$, $\delta_{s-1}=\frac{1}{2}$, and numbers $\delta_{i-1}=\eta(\kappa',\delta_i)(2\leq i\leq s-1)$ are defined in the way of Lemma \ref{sum} in successive, and $R_1=A_1$, $R_{i+1}=S_{\delta_{i}}(A_{i+1},R_i)(1\leq i\leq s-2)$. Besides, let $\min\set{2^{2-s}\delta_1,2^{3-s}\delta_2,\dots,\delta_{s-1}}\leq\eta\leq1$ be a number. Then for any arithmetic progression $P\subset[N]$ with $|P|\geq\eta N$, we must have $|P|\geq\delta_iN_i$ for all $1\leq i\leq s-1$. 

On noting that $0\leq g_i\leq1+\tilde\kappa/4$ for all $1\leq i\leq s$, for any progression $P\subset[N_{i-1}]$ (specifically,  assume that $P\subset[N_1]$ when $i=1$),
 \[
|A_i\cap P|\geq\sum_{n\in A_i\cap P}(g_i(n)-\tilde\kappa/4)=\sum_{n\in P}(g_i(n)-\tilde\kappa/4)-\sum_{n\in P\backslash A_i}(g_i(n)-\tilde\kappa/4)\geq\sum_{n\in P}g_i(n)-\frac{\tilde\kappa}{2}|P|.
\]
It then follows from the assumption $\sum_{n\in P}g_i(n)\geq(\alpha+\tilde\kappa )|P|$ whenever $|P|\geq\delta_{i-1}N_{i-1}$  (especially, $|P|\geq\delta_1N_1$ when $i=1$) that
\[
|A_i\cap P|\geq\bigbrac{\alpha+\frac{\tilde\kappa}{2}}|P|.
\]
We now apply Lemma \ref{sum} to the sets $A_1,A_2$ to get that for any progression $Q_2\subset[N_2]$ with $|Q_2|\geq\delta_2 N_2$,
\[
|R_2\cap Q_2|=|S_{\delta_1}(A_1,A_2)\cap Q_2|\geq\Bigbrac{\min\set{2\alpha+\tilde\kappa,1}-\kappa'}|Q_2|.
\]
The repeat application of Lemma \ref{sum} yields a set $R_{s-1}\subset A_1+\dots+A_{s-1}$ such that
\begin{align}\label{aaa}
|R_{s-1}\cap Q_{s-1}|\geq\Bigbrac{\min\set{(s-1)\alpha+\tilde\kappa,1}-(s-1)\kappa'}|Q_{s-1}|,
\end{align}
holds for any progression $Q_{s-1}\subset[N_{s-1}]$ with $|Q_{s-1}|\geq\delta_{s-1} N_{s-1}$

On the other hand, recalling that $R_{s-1}\subseteq A_1+\dots+A_{s-1}$, it can be seen from the definition of convolution that 
\begin{align*}
1_{A_1}*\dots*1_{A_s}(n_0)&\geq\sum_{a+b=n_0\atop a\in A_s,b\in R_{s-1}}1_{A_1}*\dots*1_{A_{s-1}}(b)1_{A_s}(a)\\
&\geq |A_s\cap(n_0-R_{s-1})|\min_{b\in R_{s-1}}1_{A_1}*\dots*1_{A_{s-1}}(b).
\end{align*}
Dues to $|A_s|\geq\sum_{n\in[N_{s-1}]}(g_s(n)-\tilde\kappa/4)-\frac{\tilde\kappa}{4}N_{s-1}\geq(\alpha+\frac{\tilde\kappa}{2})N_{s-1}$, we may take $Q_{s-1}=[N_{s-1}]$ in (\ref{aaa}) and use inclusion-exclusion principle to obtain that
\[
|A_s\cap(n_0-R_{s-1})|\geq|A_s|+|n_0-R_{s-1}|-|A_s\cup(n_0-R_{s-1})|.
\]
Thus, on noting that $A_s\cup(n_0-R_{s-1})$ is a subset of $[n_0]$,
\begin{align*}
|A_s\cap(n_0-R_{s-1})|&\geq|A_s|+|R_{s-1}|-n_0\\
&\geq\Bigbrac{\alpha+\frac{\tilde\kappa}{2}}N_{s-1}+\Bigbrac{\min\set{(s-1)\alpha+\tilde\kappa,1}-(s-1)\kappa'}N_{s-1}-n_0\\
&\geq\min\set{s\alpha+\tilde\kappa,1+\frac{\tilde\kappa}{2}}N_{s-1}-\frac{\tilde\kappa}{4}N_{s-1}-n_0.	
\end{align*}
When $s\alpha\geq1$,
\[
|A_s\cap(n_0-R_{s-1})|\geq\Bigbrac{1+\frac{\tilde\kappa}{4}}N_{s-1}-n_0.
\]
On noting that $n_0\leq N_{s-1}+2^{s-1}$,
\[
|A_s\cap(n_0-R_{s-1})|\geq\frac{\tilde\kappa}{4}N_{s-1}-2^{s-1}.
\]
Thus,
\begin{align}\label{3311}
1_{A_1}*\dots*1_{A_s}(n_0)\gg\tilde\kappa N_{s-1}\min_{b\in R_{s-1}}1_{A_1}*\dots*1_{A_{s-1}}(b).
\end{align}
Recalling the definition of $\eta$-popular sumset, 
\[
R_{s-1}=\set{b:1_{A_{s-1}}*1_{R_{s-2}}(b)\geq\delta_{s-2}\max\set{|A_{s-1}|,|R_{s-2}|}},
\] 
thus, it is not hard to find that when $b\in R_{s-1}$,
\begin{align*}
1_{A_1}*\dots*1_{A_{s-1}}(b)&\geq\sum_{a+a'=b\atop a\in A_{s-1},a'\in R_{s-2}}1_{A_1}*\dots*1_{A_{s-2}}(a')1_{A_{s-1}}(a)\\
&\geq\delta_{s-2}|A_{s-1}|\min_{a'\in R_{s-2}}1_{A_1}*\dots*1_{A_{s-2}}(a').	
\end{align*}
Following from the same manner one has for all $b\in R_{s-1}$
\begin{align}\label{3322}
1_{A_1}*\dots*1_{A_{s-1}}(b)\gg\prod_{1\leq i\leq s-2}\delta_i|A_{i+1}|\gg_{\tilde\kappa} \prod_{1\leq i\leq s-2}N_i.
\end{align}
We now replace (\ref{3322}) into (\ref{3311}), and from the assumptions  $n_0\sim N/2$ and $N_i\gg n_0$ to get that
\[
1_{A_1}*\dots*1_{A_s}(n_0)\gg_{s,\tilde\kappa}N^{s-1}.
\]

\end{document}